\documentclass[11pt]{article}
\usepackage[latin1]{inputenc}
\usepackage{amsmath,amsthm,amssymb}
\usepackage{amsfonts}
\usepackage{amsmath,amsthm,amssymb,amscd}
\usepackage{latexsym}
\usepackage{graphicx}
\usepackage{mathrsfs}
\usepackage{graphicx, color}
\usepackage{amscd}
\usepackage{amsmath}
\usepackage{amsfonts}
\usepackage{amssymb}
\usepackage{mathrsfs}

\makeatletter \@addtoreset{equation}{section} \makeatother

\newtheorem{theorem}{Theorem}[section]
\newtheorem{definition}{Definition}[section]

\newtheorem{remark} {Remark}[section]

\input{tcilatex}

\begin{document}

\title{\textbf{The Keller-Osserman type conditions for the study of
existence of entire radial solutions to a semilinear elliptic system}}
\author{ Dragos-Patru Covei$^{a}$ \\
%EndAName
{\footnotesize $^{a}$ Department of Applied Mathematics, }\\
{\footnotesize The Bucharest University of Economic Studies, }\\
{\footnotesize Piata Romana, 1st district, postal code: 010374, postal
office: 22, Romania}}
\date{ }
\maketitle

\begin{abstract}
{\ }We are concerned here with functions $\left( u_{1},u_{2}\right) $ which
are of class $C^{2}$ and satisfy the systems of the form%
\begin{equation*}
\left\{ 
\begin{array}{l}
\Delta u_{1}=p_{1}\left( \left\vert x\right\vert \right) f_{1}\left(
u_{2}\right) \text{ in }\mathbb{R}^{N}\text{, } \\ 
\Delta u_{2}=p_{2}\left( \left\vert x\right\vert \right) f_{2}\left(
u_{1}\right) \text{ in }\mathbb{R}^{N}\text{,}%
\end{array}%
\right. 
\end{equation*}%
where $p_{1}$, $f_{1}$, $p_{2}$ and $f_{2}$ are continuous functions
satisfying certain new conditions. Our considerations center about the
behavior of $\left( u_{1}\text{,}u_{2}\right) $ as $\left\vert x\right\vert
\rightarrow \infty $.
\end{abstract}

%{\section*{{\bf 1. Introduction }}} \setcounter{section}{1}
%\setcounter{equation}{0}



\section{Introduction}

Entire large and bounded solutions of semilinear elliptic systems have been
received an increased interest in the past several decades. In this article
we analyze semilinear elliptic systems of the type 
\begin{equation}
\left\{ 
\begin{array}{l}
\Delta u_{1}=p_{1}\left( \left\vert x\right\vert \right) f_{1}\left(
u_{2}\right) \text{ } \\ 
\Delta u_{2}=p_{2}\left( \left\vert x\right\vert \right) f_{2}\left(
u_{1}\right) \text{ }%
\end{array}%
\right. \text{ for }x\in \mathbb{R}^{N}\text{ (}N\geq 3\text{),}  \label{11}
\end{equation}%
in which the functions $p_{1}$, $p_{2}$, $f_{1}$ and $f_{2}$ takes various
forms, which are mentioned later on. Such problems are referred in the
literature as the Bieberbach and Rademacher problems type. The system (\ref%
{11}) will be studied under three different types of boundary conditions:

\begin{itemize}
\item \textbf{Finite Case:} Both components $(u_{1},u_{2})$ are bounded,
that is, 
\begin{equation}
\left\{ 
\begin{array}{l}
\lim_{\left\vert x\right\vert \rightarrow \infty }u_{1}\left( \left\vert
x\right\vert \right) <\infty , \\ 
\lim_{\left\vert x\right\vert \rightarrow \infty }u_{2}\left( \left\vert
x\right\vert \right) <\infty .%
\end{array}%
\right.  \label{F}
\end{equation}

\item \textbf{Infinite Case:} Both components $(u_{1},u_{2})$ are large,
that is, 
\begin{equation}
\left\{ 
\begin{array}{l}
\lim_{\left\vert x\right\vert \rightarrow \infty }u_{1}\left( \left\vert
x\right\vert \right) =\infty , \\ 
\lim_{\left\vert x\right\vert \rightarrow \infty }u_{2}\left( \left\vert
x\right\vert \right) =\infty .%
\end{array}%
\right.  \label{I}
\end{equation}

\item \textbf{Semifinite Case:} One of the components bounded while the
other is large, that is, 
\begin{equation}
\left\{ 
\begin{array}{l}
\lim_{\left\vert x\right\vert \rightarrow \infty }u_{1}\left( \left\vert
x\right\vert \right) <\infty , \\ 
\lim_{\left\vert x\right\vert \rightarrow \infty }u_{2}\left( \left\vert
x\right\vert \right) =\infty ,%
\end{array}%
\right.  \label{SF1}
\end{equation}%
or%
\begin{equation}
\left\{ 
\begin{array}{l}
\lim_{\left\vert x\right\vert \rightarrow \infty }u_{1}\left( \left\vert
x\right\vert \right) =\infty , \\ 
\lim_{\left\vert x\right\vert \rightarrow \infty }u_{2}\left( \left\vert
x\right\vert \right) <\infty .%
\end{array}%
\right.  \label{SF2}
\end{equation}
\end{itemize}

\begin{definition}
A solution $(u_{1},u_{2})\in C^{2}(\left[ 0,\infty \right) )\times C^{2}(%
\left[ 0,\infty \right) )$ of the system (\ref{11}) is called an \textit{%
entire} \textit{bounded solution} if the condition (\ref{F}) holds; an 
\textit{entire} \textit{large solution} if the condition (\ref{I}) holds; a 
\textit{semifinite entire large solution} when (\ref{SF1}) or (\ref{SF2})
holds.
\end{definition}

The study of existence of large solutions for semilinear elliptic systems of
the form (\ref{11}) goes back to the pioneering papers by Keller \cite{K}
and Osserman \cite{O}. In 1957 Osserman \cite{O} proved that, for a given
positive, continuous and nondecreasing function $f$, the semilinear elliptic
partial differential inequality 
\begin{equation}
\Delta u\geq a\left( x\right) f\left( u\right) \text{ in }\mathbb{R}^{N},
\label{os}
\end{equation}%
with $a\left( x\right) =1$, possesses an entire large solution $u:\mathbb{R}%
^{N}\rightarrow \mathbb{R}$ if and only if 
\begin{equation}
\int_{1}^{\infty }\left( \int_{0}^{t}f\left( s\right) ds\right)
^{-1/2}dt=+\infty ,  \label{ko}
\end{equation}%
(assumption known today as the Keller-Osserman condition).

Such problems are drawn by the mathematical modelling of many natural
phenomena related to steady-state reaction-diffusion, subsonic fluid flows,
electrostatic potential in a shiny metallic body inside or subsonic motion
of a gas, automorphic functions theory, geometry and control theory (see,
for example, L. Bieberbach \cite{LB}, Diaz \cite{DZ}, Keller \cite{KG},
Lasry and Lions \cite{LAS}, Matero \cite{JM}, Marinescu and Varsan \cite{M1}%
, Iftimie-Marinescu and Varsan \cite{M2},\textit{\ }and Rademacher\textit{\ }%
\cite{R}\textit{\ }for a more detailed discussion). For example, reading the
work of Lasry and Lions \cite{LAS}, we can observe that such problems arise
in stochastic control theory. The controls are to be designed so that the
state of the system is constrained to some region. Finding optimal controls
is then shown to be equivalent to finding large solutions for a second order
nonlinear elliptic partial differential equation.

Even if these problems are treated directly or indirectly, in many papers
from the specialized literature, they have not been completely clarified
yet. This does not surprise us, because practical applications mainly reveal
new horizons, complexity and aspects that allow new theoretical approaches.

Our objective in the present research, in short, is to complete and to find
new ideas to treat the principal results of the author \cite{CDF,CD2,CD3},
Goyal \cite{VG,VGP}, Magliaro-Mari-Mastrolia and Rigoli \cite{MAR},
Lieberman \cite{GL}, Nehari \cite{NH}, Redheffer \cite{RR}, Rhee \cite{RH},
Reichel-Walter \cite{RW} and other associated works.

More exactly, the present work is related to the first and last of the
preceding questions:

1. The first one is the problem of existence of solution to (\ref{11}) that
satisfies (\ref{F}), (\ref{I}), (\ref{SF1}) or (\ref{SF2})

2. The second one is to give a necessary and a sufficient condition for a
positive radial solution of (\ref{11}) to be entire large.

However, there is no results for systems (\ref{11}), where $f_{1}$, $f_{2}$
satisfy a condition of the form (\ref{ko}). Other purpose of this paper is
to fill this gap. To be more precise, we consider $a$, $b\in \left( 0,\infty
\right) $ arbitrary and we assume that the variable weights functions $%
p_{1},p_{2}$ and the nonlinearities $f_{1}$, $f_{2}$ satisfy:

(P1)\quad $p_{1},p_{2}:\left[ 0,\infty \right) \rightarrow \left[ 0,\infty
\right) $ are spherically symmetric continuous functions (i.e.,\textit{\ }$%
p_{1}\left( x\right) =p_{1}\left( \left\vert x\right\vert \right) $ and $%
p_{2}\left( x\right) =p_{2}\left( \left\vert x\right\vert \right) $)\textit{;%
}

(C1)\quad $f_{1}$, $f_{2}:\left[ 0,\infty \right) \rightarrow \left[
0,\infty \right) $ are continuous, non-decreasing, $f_{1}\left( 0\right) =$ $%
f_{2}\left( 0\right) =0$ and $f_{1}\left( s\right) >0$, $f_{2}\left(
s\right) >0$ for all $s>0$;

(C2)\quad there exist positive constants $\overline{c}_{1},\overline{c}_{2}$%
, the continuous and increasing functions $h_{1}$, $h_{2}$, $\omega _{1}$, $%
\omega _{2}:\left[ 0,\infty \right) \rightarrow \left[ 0,\infty \right) $
such that%
\begin{eqnarray}
f_{1}\left( t_{1}\cdot w_{1}\right) &\leq &\overline{c}_{1}h_{1}\left(
t_{1}\right) \cdot \omega _{1}\left( w_{1}\right) \text{ }\forall \text{ }%
w_{1}\geq 1\text{ and }\forall \text{ }t_{1}\geq M_{1}\cdot f_{2}\left(
a\right) ,  \label{c21} \\
f_{2}\left( t_{2}\cdot w_{2}\right) &\leq &\overline{c}_{2}h_{2}\left(
t_{2}\right) \cdot \omega _{2}\left( w_{2}\right) \text{ }\forall \text{ }%
w_{2}\geq 1\text{ and }\forall \text{ }t_{2}\geq M_{2}\cdot f_{1}\left(
b\right) ,  \label{c22}
\end{eqnarray}%
where $M_{1}$ and $M_{2}$ are such that 
\begin{equation*}
M_{1}=\left\{ 
\begin{array}{lll}
\frac{b}{f_{2}\left( a\right) } & if & b>f_{2}\left( a\right) , \\ 
1 & if & \text{ }b\leq f_{2}\left( a\right) ,%
\end{array}%
\right. \text{ \ and\ \ }M_{2}=\left\{ 
\begin{array}{lll}
\frac{a}{f_{1}\left( b\right) } & if & \text{ }a>f_{1}\left( b\right) , \\ 
1 & if & \text{ }a\leq f_{1}\left( b\right) .%
\end{array}%
\right.
\end{equation*}%
To facilitate the presentation of the results we introduce some notations:%
\begin{eqnarray*}
G_{1}\left( z\right) &=&\int_{0}^{z}s^{N-1}p_{2}\left( s\right) ds\text{ } \\
G_{2}\left( z\right) &=&\int_{0}^{z}s^{N-1}p_{1}\left( s\right) ds \\
P_{1}\left( r\right) &=&\int_{0}^{r}y^{1-N}\int_{0}^{y}t^{N-1}p_{1}\left(
t\right) f_{1}\left( b+f_{2}\left( a\right) \int_{0}^{t}z^{1-N}G_{1}\left(
z\right) dz\right) dtdy, \\
Q_{1}\left( r\right) &=&\int_{0}^{r}y^{1-N}\int_{0}^{y}t^{N-1}p_{2}\left(
t\right) f_{2}\left( a+f_{1}\left( b\right) \int_{0}^{t}z^{1-N}G_{2}\left(
z\right) dz\right) dtdy, \\
P_{2}\left( r\right) &=&\int_{0}^{r}z^{1+\varepsilon }p_{1}\left( z\right)
\omega _{1}\left( 1+\int_{0}^{z}\frac{1}{t^{N-1}}\int_{0}^{t}s^{N-1}p_{2}%
\left( s\right) dsdt\right) dz, \\
Q_{2}\left( r\right) &=&\int_{0}^{r}z^{1+\varepsilon }p_{2}\left( z\right)
\omega _{2}\left( 1+\int_{0}^{z}t^{1-N}\int_{0}^{t}s^{N-1}p_{1}\left(
s\right) dsdt\right) dz, \\
P_{3}\left( r\right) &=&\int_{0}^{r}\sqrt{2\phi _{1}\left( z\right) \omega
_{1}\left( 1+\int_{0}^{z}t^{1-N}\int_{0}^{t}s^{N-1}p_{2}\left( s\right)
dsdt\right) }dz\text{ where }\phi _{1}\left( z\right) =\underset{0\leq t\leq
z}{\max }p_{1}\left( t\right) , \\
Q_{3}\left( r\right) &=&\int_{0}^{r}\sqrt{2\phi _{2}\left( z\right) \omega
_{2}\left( 1+\int_{0}^{z}t^{1-N}\int_{0}^{t}s^{N-1}p_{1}\left( s\right)
dsdt\right) }dz\text{ where }\phi _{2}\left( z\right) =\underset{0\leq t\leq
z}{\max }p_{2}\left( t\right) , \\
H_{1}\left( r\right) &=&\int_{a}^{r}\frac{1}{\sqrt{\int_{0}^{s}h_{1}\left(
M_{1}f_{2}\left( t\right) \right) dt}}ds\text{, }H_{2}\left( r\right)
=\int_{b}^{r}\frac{1}{\sqrt{\int_{0}^{s}h_{2}\left( M_{2}f_{1}\left(
t\right) \right) dt}}ds, \\
P_{i}\left( \infty \right) &=&\lim_{r\rightarrow \infty }P_{i}\left(
r\right) ,\text{ }Q_{i}\left( \infty \right) =\lim_{r\rightarrow \infty
}Q_{i}\left( r\right) \text{, }H_{i}\left( \infty \right)
=\lim_{r\rightarrow \infty }H_{i}\left( r\right) \text{ for }i=1,2.
\end{eqnarray*}

Our main result are the following:

\begin{theorem}
\label{th1}\textit{Assume that }$H_{1}\left( \infty \right) =H_{2}\left(
\infty \right) =\infty $\textit{. If }$p_{1}$, $p_{2}$,\ $f_{1}$, $f_{2}$ 
\textit{\ satisfy \textrm{(P1)}},\textit{\ \textrm{(C1)} and \textrm{(C2)},
then the problem \textrm{(\ref{11})} has a nonnegative entire radial
solution }$\left( u_{1},u_{2}\right) $ with central value in $\left(
a,b\right) $ (i.e. $\left( u_{1},u_{2}\right) =\left( a,b\right) $)\textit{. 
}Moreover,

i.)\quad if in addition\textit{\ }$r^{2N-2}p_{1}\left( r\right) $\textit{, }$%
r^{2N-2}p_{2}\left( r\right) $\textit{\ are nondecreasing for large }$r$%
\textit{\ and} $p_{1}$, $p_{2}$ satisfy 
\begin{equation*}
P_{2}\left( \infty \right) <\infty \text{ and }Q_{2}\left( \infty \right)
<\infty ,
\end{equation*}%
then for any nonnegative radial solution $\left( u_{1},u_{2}\right) $ of 
\textrm{(\ref{11}) }with central value in $\left( a,b\right) $ we have (\ref%
{F});

ii.)\quad if $p_{1}$ and $p_{2}$ satisfy 
\begin{equation}
P_{1}\left( \infty \right) =\text{ }Q_{1}\left( \infty \right) =\infty ,
\label{12}
\end{equation}%
and $\left( u_{1},u_{2}\right) $ is any nonnegative radial solution of 
\textrm{(\ref{11}) }with central value in $\left( a,b\right) $ then (\ref{I}%
) holds;

iii.)\quad if in addition\textit{\ }$r^{2N-2}p_{1}\left( r\right) $\textit{\
is nondecreasing for large }$r$\textit{\ and} $p_{1}$, $p_{2}$ satisfy 
\begin{equation}
P_{2}\left( \infty \right) <\infty \text{ and }Q_{1}\left( \infty \right)
=\infty ,\text{ }  \label{123}
\end{equation}%
then for any nonnegative radial solution $\left( u_{1},u_{2}\right) $ of 
\textrm{(\ref{11}) }with central value in $\left( a,b\right) $ we have (\ref%
{SF1});

iv.)\quad if in addition\textit{\ }$r^{2N-2}p_{2}\left( r\right) $\textit{\
is nondecreasing for large }$r$\textit{\ and} $p_{1}$, $p_{2}$ satisfy 
\begin{equation}
P_{1}\left( \infty \right) =\infty \text{ and }Q_{2}\left( \infty \right)
<\infty ,  \label{122}
\end{equation}%
then for any nonnegative radial solution $\left( u_{1},u_{2}\right) $ of 
\textrm{(\ref{11}) }with central value in $\left( a,b\right) $ we have (\ref%
{SF2});

v.)\quad if \textrm{(\ref{11})} has a nonnegative entire large solution $%
\left( u_{1},u_{2}\right) $ with central value in $\left( a,b\right) $ and $%
r^{2N-2}p_{1}\left( r\right) $\textit{, }$r^{2N-2}p_{2}\left( r\right) $%
\textit{\ are nondecreasing for large }$r$, then $p_{1}$ and $p_{2}$ satisfy%
\begin{equation}
P_{2}\left( \infty \right) =Q_{2}\left( \infty \right) =\infty ,  \label{13}
\end{equation}%
for every $\varepsilon >0$.
\end{theorem}

\begin{theorem}
\label{th2}Assume that the hypotheses \textrm{(P1),} \textrm{(C1) }and%
\textrm{\ (C2)} are satisfied. \ The following hold:

i.)\quad\ If $P_{3}\left( \infty \right) <H_{1}\left( \infty \right) <\infty 
$ and $Q_{3}\left( \infty \right) <H_{2}\left( \infty \right) <\infty $ then
the system (\ref{11}) has one positive bounded radial solution $\left(
u_{1},u_{2}\right) \in C^{2}\left( \left[ 0,\infty \right) \right) \times
C^{2}\left( \left[ 0,\infty \right) \right) $, with central value in $\left(
a,b\right) ,$ such that%
\begin{equation*}
\left\{ 
\begin{array}{c}
a+P_{1}\left( r\right) \leq u_{1}\left( r\right) \leq H_{1}^{-1}\left( 
\overline{c}_{1}^{1/2}P_{3}\left( r\right) \right) , \\ 
b+Q_{1}\left( r\right) \leq u_{2}\left( r\right) \leq H_{2}^{-1}\left( 
\overline{c}_{2}^{1/2}Q_{3}\left( r\right) \right) .%
\end{array}%
\right. \text{ }
\end{equation*}

ii.)\quad If $H_{1}\left( \infty \right) =\infty $, $P_{1}\left( \infty
\right) =\infty $ and $Q_{3}\left( \infty \right) <H_{2}\left( \infty
\right) <\infty $ then the system (\ref{11}) has one positive radial
solution 
\begin{equation*}
\left( u_{1},u_{2}\right) \in C^{2}\left( \left[ 0,\infty \right) \right)
\times C^{2}\left( \left[ 0,\infty \right) \right) ,
\end{equation*}%
with central value in $\left( a,b\right) $, such that (\ref{SF2}) holds;

iii.)\quad If $P_{3}\left( \infty \right) <H_{1}\left( \infty \right)
<\infty $ and $H_{2}\left( \infty \right) =\infty $, $Q_{1}\left( \infty
\right) =\infty $ then the system (\ref{11}) has one positive radial
solution 
\begin{equation*}
\left( u_{1},u_{2}\right) \in C^{2}\left( \left[ 0,\infty \right) \right)
\times C^{2}\left( \left[ 0,\infty \right) \right) ,
\end{equation*}%
with central value in $\left( a,b\right) $, such that (\ref{SF1}) holds;

iv)\quad If \textit{\ }$r^{2N-2}p_{1}\left( r\right) $\textit{\ is
nondecreasing for large }$r$, $H_{1}\left( \infty \right) =\infty $, $%
P_{2}\left( \infty \right) <\infty $ and $Q_{3}\left( \infty \right)
<H_{2}\left( \infty \right) <\infty $ then the system (\ref{11}) has one
positive radial solution 
\begin{equation*}
\left( u_{1},u_{2}\right) \in C^{2}\left( \left[ 0,\infty \right) \right)
\times C^{2}\left( \left[ 0,\infty \right) \right) ,
\end{equation*}%
with central value in $\left( a,b\right) $, such that (\ref{F}) holds;

v)\quad If \textit{\ }$r^{2N-2}p_{2}\left( r\right) $\textit{\ is
nondecreasing for large }$r$, $P_{3}\left( \infty \right) <H_{1}\left(
\infty \right) <\infty $ and $H_{2}\left( \infty \right) =\infty $, $%
Q_{2}\left( \infty \right) <\infty $ then the system (\ref{11}) has one
positive radial solution 
\begin{equation*}
\left( u_{1},u_{2}\right) \in C^{2}\left( \left[ 0,\infty \right) \right)
\times C^{2}\left( \left[ 0,\infty \right) \right) ,
\end{equation*}%
with central value in $\left( a,b\right) $, such that (\ref{F}) holds;
\end{theorem}

\begin{remark}
Our assumption (C2) is further discussed in the famous book of
Krasnosel'skii and Rutickii \cite{KR} (see also Gustavsson, Maligranda and
Peetre \cite{GU}).
\end{remark}

\begin{remark}
A simple example of nonlinearities $f_{1}$ and $f_{2}$ satisfying the
assumptions in Theorem \ref{th1} and Theorem \ref{th2} are%
\begin{equation*}
f_{1}\left( u_{2}\right) =h_{1}\left( u_{2}\right) =\omega _{1}\left(
u_{2}\right) =u_{2}^{\alpha }\text{ and }f_{2}\left( u_{1}\right)
=h_{2}\left( u_{1}\right) =\omega _{1}\left( u_{1}\right) =u_{1}^{\beta }
\end{equation*}%
where $\alpha ,\beta \in \mathbb{R}$. The results in Theorem \ref{th1} work
with $\alpha \cdot \beta \leq 1$ and $\overline{c}_{1}=$ $\overline{c}%
_{2}=1. $ On the other hand the results in Theorem \ref{th2} are proved for $%
\alpha \cdot \beta >1$ and $\overline{c}_{1}=$ $\overline{c}_{2}=1$.
\end{remark}

\section{Proofs of the Theorems}

We prove the existence of a solution $\left( u_{1},u_{2}\right) $ for the
system%
\begin{equation}
\left\{ 
\begin{array}{l}
\Delta u_{1}\left( r\right) =p_{1}\left( r\right) f_{1}\left( u_{2}\left(
r\right) \right) \text{ for }r:=\left\vert x\right\vert \text{,} \\ 
\Delta u_{2}\left( r\right) =p_{2}\left( r\right) f_{2}\left( u_{1}\left(
r\right) \right) \text{ for }r:=\left\vert x\right\vert \text{.}%
\end{array}%
\right.  \label{6}
\end{equation}%
In the radial setting, the system (\ref{6}) becomes a system of differential
equations of the form 
\begin{equation}
\left\{ 
\begin{array}{l}
\left( r^{N-1}u_{1}^{\prime }\left( r\right) \right) ^{\prime
}=r^{N-1}p_{1}\left( r\right) f_{1}\left( u_{2}\left( r\right) \right) , \\ 
\left( r^{N-1}u_{2}^{\prime }\left( r\right) \right) ^{\prime
}=r^{N-1}p_{2}\left( r\right) f_{2}\left( u_{1}\left( r\right) \right) \text{
}.%
\end{array}%
\right.  \label{77}
\end{equation}%
Then, a radial solution of (\ref{77}) is any solution $\left(
u_{1},u_{2}\right) $ for the integral equations%
\begin{equation*}
\left\{ 
\begin{array}{l}
u_{1}\left( r\right) =a+\int_{0}^{r}t^{1-N}\int_{0}^{t}s^{N-1}p_{1}\left(
s\right) f_{1}\left( u_{2}\left( s\right) \right) dsdt,\text{ } \\ 
u_{2}\left( r\right) =b+\int_{0}^{r}t^{1-N}\int_{0}^{t}s^{N-1}p_{2}\left(
s\right) f_{2}\left( u_{1}\left( s\right) \right) \text{ }dsdt.%
\end{array}%
\right.
\end{equation*}%
To establish a solution to this system, we use successive approximation.
Define sequences $\left\{ u_{1}^{k}\right\} ^{k\geq 1}$ and $\left\{
u_{2}^{k}\right\} ^{k\geq 1}$ on $\left[ 0,\infty \right) $ by%
\begin{equation*}
\left\{ 
\begin{array}{l}
u_{1}^{0}=a\text{ and }u_{2}^{0}=b,\text{ }r\geq 0, \\ 
u_{1}^{k}\left( r\right)
=a+\int_{0}^{r}t^{1-N}\int_{0}^{t}s^{N-1}p_{1}\left( s\right) f_{1}\left(
u_{2}^{k-1}\left( s\right) \right) dsdt, \\ 
u_{2}^{k}\left( r\right)
=b+\int_{0}^{r}t^{1-N}\int_{0}^{t}s^{N-1}p_{2}\left( s\right) f_{2}\left(
u_{1}^{k-1}\left( s\right) \right) dsdt.%
\end{array}%
\right.
\end{equation*}%
We remark that, for all $r\geq 0$ and $k\in N$ 
\begin{equation*}
u_{1}^{k}\left( r\right) \geq a\text{ and }u_{2}^{k}\left( r\right) \geq b%
\text{.}
\end{equation*}%
Moreover, proceeding by mathematical induction we conclude that $\left\{
u_{1}^{k}\right\} ^{k\geq 1}$ and $\left\{ u_{2}^{k}\right\} ^{k\geq 1}$ are
non-decreasing sequences on $\left[ 0,\infty \right) $. We will next prove
the \textquotedblright upper bounds\textquotedblright . To do this, we note
that $\left\{ u_{1}^{k}\right\} ^{k\geq 1}$ and $\left\{ u_{2}^{k}\right\}
^{k\geq 1}$ satisfy%
\begin{equation*}
\left\{ 
\begin{array}{l}
\left[ r^{N-1}\left( u_{1}^{k}\left( r\right) \right) ^{\prime }\right]
^{\prime }=r^{N-1}p_{1}\left( r\right) f_{1}\left( u_{2}^{k-1}\left(
r\right) \right) ,\text{ } \\ 
\left[ r^{N-1}\left( u_{2}^{k}\left( r\right) \right) ^{\prime }\right]
^{\prime }=r^{N-1}p_{2}\left( r\right) f_{2}\left( u_{1}^{k-1}\left(
r\right) \right) .%
\end{array}%
\right.
\end{equation*}%
Using the monotonicity of $\left\{ u_{1}^{k}\right\} ^{k\geq 1}$ and $%
\left\{ u_{2}^{k}\right\} ^{k\geq 1}$ we find the inequalities%
\begin{equation}
\left\{ 
\begin{array}{l}
\left[ r^{N-1}\left( u_{1}^{k}\left( r\right) \right) ^{\prime }\right]
^{\prime }=r^{N-1}p_{1}\left( r\right) f_{1}\left( u_{2}^{k-1}\left(
r\right) \right) \leq r^{N-1}p_{1}\left( r\right) f_{1}\left(
u_{2}^{k}\left( r\right) \right) , \\ 
\left[ r^{N-1}\left( u_{2}^{k}\left( r\right) \right) ^{\prime }\right]
^{\prime }=r^{N-1}p_{2}\left( r\right) f_{2}\left( u_{1}^{k-1}\left(
r\right) \right) \leq r^{N-1}p_{2}\left( r\right) f_{2}\left(
u_{1}^{k}\left( r\right) \right) .%
\end{array}%
\right.  \label{8}
\end{equation}%
Then, going back to the previous computation we have%
\begin{eqnarray}
\left[ r^{N-1}\left( u_{1}^{k}\left( r\right) \right) ^{\prime }\right]
^{\prime } &=&r^{N-1}p_{1}\left( r\right) f_{1}\left( u_{2}^{k-1}\left(
r\right) \right)  \notag \\
&\leq &r^{N-1}p_{1}\left( r\right) f_{1}\left( u_{2}^{k}\left( r\right)
\right)  \notag \\
&\leq &r^{N-1}p_{1}\left( r\right) f_{1}\left(
b+\int_{0}^{r}t^{1-N}\int_{0}^{t}s^{N-1}p_{2}\left( s\right) f_{2}\left(
u_{1}^{k-1}\left( s\right) \right) dsdt\right)  \notag \\
&\leq &r^{N-1}p_{1}\left( r\right) f_{1}\left(
b+\int_{0}^{r}t^{1-N}\int_{0}^{t}s^{N-1}p_{2}\left( s\right) f_{2}\left(
u_{1}^{k}\left( s\right) \right) dsdt\right)  \notag \\
&\leq &r^{N-1}p_{1}\left( r\right) f_{1}\left( b+f_{2}\left( u_{1}^{k}\left(
r\right) \right) \int_{0}^{r}t^{1-N}\int_{0}^{t}s^{N-1}p_{2}\left( s\right)
dsdt\right)  \label{n8} \\
&\leq &r^{N-1}p_{1}\left( r\right) f_{1}\left( f_{2}\left( u_{1}^{k}\left(
r\right) \right) \left( \frac{b}{f_{2}\left( u_{1}^{k}\left( r\right)
\right) }+\int_{0}^{r}t^{1-N}\int_{0}^{t}s^{N-1}p_{2}\left( s\right)
dsdt\right) \right)  \notag \\
&\leq &r^{N-1}p_{1}\left( r\right) f_{1}\left( f_{2}\left( u_{1}^{k}\left(
r\right) \right) \left( \frac{b}{f\left( a\right) }+\int_{0}^{r}t^{1-N}%
\int_{0}^{t}s^{N-1}p_{2}\left( s\right) dsdt\right) \right)  \notag \\
&\leq &r^{N-1}p_{1}\left( r\right) f_{1}\left( M_{1}f_{2}\left(
u_{1}^{k}\left( r\right) \right) \left(
1+\int_{0}^{r}t^{1-N}\int_{0}^{t}s^{N-1}p_{2}\left( s\right) dsdt\right)
\right)  \notag \\
&\leq &r^{N-1}p_{1}\left( r\right) \overline{c}_{1}h_{1}\left(
M_{1}f_{2}\left( u_{1}^{k}\left( r\right) \right) \right) \omega _{1}\left(
1+\int_{0}^{r}t^{1-N}\int_{0}^{t}s^{N-1}p_{2}\left( s\right) dsdt\right) , 
\notag
\end{eqnarray}%
and%
\begin{eqnarray}
\left[ r^{N-1}\left( u_{2}^{k}\left( r\right) \right) ^{\prime }\right]
^{\prime } &=&r^{N-1}p_{2}\left( r\right) f_{2}\left( u_{1}^{k-1}\left(
r\right) \right)  \notag \\
&\leq &r^{N-1}p_{2}\left( r\right) f_{2}\left( M_{2}f_{1}\left(
u_{2}^{k}\left( r\right) \right) \left(
1+\int_{0}^{r}t^{1-N}\int_{0}^{t}s^{N-1}p_{1}\left( s\right) dsdt\right)
\right)  \label{nn8} \\
&\leq &r^{N-1}p_{2}\left( r\right) \overline{c}_{2}h_{2}\left(
M_{2}f_{1}\left( u_{2}^{k}\left( r\right) \right) \right) \omega _{2}\left(
1+\int_{0}^{r}t^{1-N}\int_{0}^{t}s^{N-1}p_{1}\left( s\right) dsdt\right) . 
\notag
\end{eqnarray}%
By (\ref{n8}) and (\ref{nn8}), we have%
\begin{equation}
\left\{ 
\begin{array}{l}
r^{N-1}\left( u_{1}^{k}\right) ^{\prime \prime }\leq \left( N-1\right)
r^{N-2}\left( u_{1}^{k}\right) ^{\prime }+r^{N-1}\left( u_{1}^{k}\right)
^{\prime \prime }=\left[ r^{N-1}\left( u_{1}^{k}\right) ^{\prime }\right]
^{\prime } \\ 
\leq r^{N-1}p_{1}\left( r\right) \overline{c}_{1}h_{1}\left(
M_{1}f_{2}\left( u_{1}^{k}\left( r\right) \right) \right) \omega _{1}\left(
1+\int_{0}^{r}t^{1-N}\int_{0}^{t}s^{N-1}p_{2}\left( s\right) dsdt\right) ,
\\ 
r^{N-1}\left( u_{2}^{k}\right) ^{\prime \prime }\leq \left[ r^{N-1}\left(
u_{2}^{k}\right) ^{\prime }\right] ^{\prime } \\ 
\leq r^{N-1}p_{2}\left( r\right) \overline{c}_{2}h_{2}\left(
M_{2}f_{1}\left( u_{2}^{k}\left( r\right) \right) \right) \omega _{2}\left(
1+\int_{0}^{r}t^{1-N}\int_{0}^{t}s^{N-1}p_{1}\left( s\right) dsdt\right) ,%
\end{array}%
\right.  \label{88}
\end{equation}%
Multiplying the first inequality by $\left( u_{1}^{k}\right) ^{\prime }$ and
the second by $\left( u_{2}^{k}\right) ^{\prime }$, we obtain 
\begin{equation}
\left\{ 
\begin{array}{l}
\left\{ \left[ \left( u_{1}^{k}\left( r\right) \right) ^{\prime }\right]
^{2}\right\} ^{\prime }\leq 2p_{1}\left( r\right) \overline{c}%
_{1}h_{1}\left( M_{1}f_{2}\left( u_{1}^{k}\left( r\right) \right) \right)
\left( u_{1}^{k}\left( r\right) \right) ^{\prime }\omega _{1}\left(
1+\int_{0}^{r}t^{1-N}\int_{0}^{t}s^{N-1}p_{2}\left( s\right) dsdt\right) ,
\\ 
\left\{ \left[ \left( u_{2}^{k}\left( r\right) \right) ^{\prime }\right]
^{2}\right\} ^{\prime }\leq 2p_{2}\left( r\right) \overline{c}%
_{2}h_{2}\left( M_{2}f_{1}\left( u_{2}^{k}\left( r\right) \right) \right)
\left( u_{2}^{k}\left( r\right) \right) ^{\prime }\omega _{2}\left(
1+\int_{0}^{r}t^{1-N}\int_{0}^{t}s^{N-1}p_{1}\left( s\right) dsdt\right) .%
\end{array}%
\right.  \label{for}
\end{equation}%
Integrating in (\ref{for}) from $0$ to $r$ we also have%
\begin{equation}
\left\{ 
\begin{array}{l}
\left[ \left( u_{1}^{k}\left( r\right) \right) ^{\prime }\right] ^{2}\leq
\int_{0}^{r}2p_{1}\left( z\right) \overline{c}_{1}h_{1}\left(
M_{1}f_{2}\left( u_{1}^{k}\left( z\right) \right) \right) \left(
u_{1}^{k}\left( z\right) \right) ^{\prime }\omega _{1}\left(
1+\int_{0}^{z}t^{1-N}\int_{0}^{t}s^{N-1}p_{2}\left( s\right) dsdt\right) dz,
\\ 
\left[ \left( u_{2}^{k}\left( r\right) \right) ^{\prime }\right] ^{2}\leq
\int_{0}^{r}2p_{2}\left( z\right) \overline{c}_{2}h_{2}\left(
M_{2}f_{1}\left( u_{2}^{k}\left( z\right) \right) \right) \left(
u_{2}^{k}\left( z\right) \right) ^{\prime }\omega _{2}\left(
1+\int_{0}^{z}t^{1-N}\int_{0}^{t}s^{N-1}p_{1}\left( s\right) dsdt\right) dz.%
\end{array}%
\right.  \label{ineqr}
\end{equation}%
Set now 
\begin{eqnarray}
\phi _{1}\left( r\right) &=&\max \left\{ p_{1}\left( z\right) \left\vert
0\leq z\leq r\right. \right\} \text{,}  \label{mom2} \\
\phi _{2}\left( r\right) &=&\max \left\{ p_{2}\left( z\right) \left\vert
0\leq z\leq r\right. \right\} \text{.}  \notag
\end{eqnarray}%
Thanks to the definition of $\phi _{1}\left( r\right) $ and $\phi _{2}\left(
r\right) $ we get from the inequalities (\ref{ineqr}) that%
\begin{equation}
\left\{ 
\begin{array}{l}
\left[ \left( u_{1}^{k}\left( r\right) \right) ^{\prime }\right] ^{2}\leq 2%
\overline{c}_{1}\phi _{1}\left( r\right) \omega _{1}\left(
1+\int_{0}^{r}t^{1-N}\int_{0}^{t}s^{N-1}p_{2}\left( s\right) dsdt\right)
\int_{0}^{r}h_{1}\left( M_{1}f_{2}\left( u_{1}^{k}\left( z\right) \right)
\right) \left( u_{1}^{k}\left( z\right) \right) ^{\prime }dz\text{, } \\ 
\left[ \left( u_{2}^{k}\left( r\right) \right) ^{\prime }\right] ^{2}\leq 2%
\overline{c}_{2}\phi _{2}\left( r\right) \omega _{2}\left(
1+\int_{0}^{r}t^{1-N}\int_{0}^{t}s^{N-1}p_{1}\left( s\right) dsdt\right)
\int_{0}^{r}h_{2}\left( M_{2}f_{1}\left( u_{2}^{k}\left( z\right) \right)
\right) \left( u_{2}^{k}\left( z\right) \right) ^{\prime }dz.%
\end{array}%
\right.  \label{mom}
\end{equation}%
As a consequence of (\ref{mom}), we also have%
\begin{equation}
\left\{ 
\begin{array}{l}
\left( u_{1}^{k}\left( r\right) \right) ^{\prime }\leq \sqrt{2\phi
_{1}\left( r\right) \omega _{1}\left(
1+\int_{0}^{r}t^{1-N}\int_{0}^{t}s^{N-1}p_{2}\left( s\right) dsdt\right) }%
\left( \int_{a}^{u_{1}^{k}\left( r\right) }\overline{c}_{1}h_{1}\left(
M_{1}f_{2}\left( z\right) \right) dz\right) ^{1/2}\text{, } \\ 
\left( u_{2}^{k}\left( r\right) \right) ^{\prime }\leq \sqrt{2\phi
_{2}\left( r\right) \omega _{2}\left(
1+\int_{0}^{r}t^{1-N}\int_{0}^{t}s^{N-1}p_{1}\left( s\right) dsdt\right) }%
\left( \int_{b}^{u_{2}^{k}\left( r\right) }\overline{c}_{2}h_{2}\left(
M_{2}f_{1}\left( z\right) \right) dz\right) ^{1/2},%
\end{array}%
\right.  \label{inte}
\end{equation}%
and, thus%
\begin{equation}
\left\{ 
\begin{array}{l}
\frac{\left( u_{1}^{k}\left( r\right) \right) ^{\prime }}{\left(
\int_{a}^{u_{1}^{k}\left( r\right) }h_{1}\left( M_{1}f_{2}\left( z\right)
\right) dz\right) ^{1/2}}\leq \overline{c}_{1}^{1/2}\sqrt{2\phi _{1}\left(
r\right) \omega _{1}\left(
1+\int_{0}^{r}t^{1-N}\int_{0}^{t}s^{N-1}p_{2}\left( s\right) dsdt\right) },
\\ 
\frac{\left( u_{2}^{k}\left( r\right) \right) ^{\prime }}{\left(
\int_{b}^{u_{2}^{k}\left( r\right) }h_{2}\left( M_{2}f_{1}\left( z\right)
\right) dz\right) ^{1/2}}\leq \overline{c}_{2}^{1/2}\sqrt{2\phi _{2}\left(
r\right) \omega _{2}\left(
1+\int_{0}^{r}t^{1-N}\int_{0}^{t}s^{N-1}p_{1}\left( s\right) dsdt\right) }.%
\end{array}%
\right.  \label{inte2}
\end{equation}%
Integrating (\ref{inte2}) leads to 
\begin{equation*}
\left\{ 
\begin{array}{l}
\int_{a}^{u_{1}^{k}\left( r\right) }\frac{1}{\sqrt{\int_{0}^{s}h_{1}\left(
M_{1}f_{2}\left( t\right) \right) dt}}ds\leq \overline{c}_{1}^{1/2}%
\int_{0}^{r}\sqrt{2\phi _{1}\left( z\right) \omega _{1}\left(
1+\int_{0}^{z}t^{1-N}\int_{0}^{t}s^{N-1}p_{2}\left( s\right) dsdt\right) }dz,
\\ 
\int_{b}^{u_{2}^{k}\left( r\right) }\frac{1}{\sqrt{\int_{0}^{s}h_{2}\left(
M_{2}f_{1}\left( t\right) \right) dt}}ds\leq \overline{c}_{2}^{1/2}%
\int_{0}^{r}\sqrt{2\phi _{2}\left( z\right) \omega _{2}\left(
1+\int_{0}^{z}t^{1-N}\int_{0}^{t}s^{N-1}p_{1}\left( s\right) dsdt\right) }dz,%
\end{array}%
\right.
\end{equation*}%
which is equivalent to%
\begin{equation}
\left\{ 
\begin{array}{l}
H_{1}\left( u_{1}^{k}\left( r\right) \right) \leq \overline{c}%
_{1}^{1/2}P_{3}\left( r\right) , \\ 
H_{2}\left( u_{2}^{k}\left( r\right) \right) \leq \overline{c}%
_{2}^{1/2}Q_{3}\left( r\right) .%
\end{array}%
\right.  \label{pen}
\end{equation}%
Since $H_{1}^{-1}$ and $H_{2}^{-1}$ are strictly increasing on $\left[
0,\infty \right) $, as previously discussed, we have that%
\begin{equation}
\left\{ 
\begin{array}{l}
u_{1}^{k}\left( r\right) \leq H_{1}^{-1}\left( \overline{c}%
_{1}^{1/2}P_{3}\left( r\right) \right) , \\ 
u_{2}^{k}\left( r\right) \leq H_{2}^{-1}\left( \overline{c}%
_{2}^{1/2}Q_{3}\left( r\right) \right) .%
\end{array}%
\right.  \label{flori}
\end{equation}%
These inequalities are independent of $k$.

\subparagraph{\textbf{Proof of Theorem \protect\ref{th1} completed}}

Combining (\ref{pen}) with 
\begin{equation*}
H_{1}\left( \infty \right) =H_{2}\left( \infty \right) =\infty
\end{equation*}%
yields that the sequences $\left\{ u_{1}^{k}\right\} ^{k\geq 1}$ and $%
\left\{ u_{2}^{k}\right\} ^{k\geq 1}$ are bounded and equicontinuous on $%
\left[ 0,c_{0}\right] $ for arbitrary $c_{0}>0$. Possibly after passing to a
subsequence, we may assume that $\left\{ \left( u_{1}^{k},u_{2}^{k}\right)
\right\} ^{^{k\geq 1}}$ converges uniformly to $\left( u_{1},u_{2}\right) $
on $\left[ 0,c_{0}\right] $. At the end of this process, we conclude by the
arbitrariness of $c_{0}>0$, that $\left( u_{1},u_{2}\right) $ is a positive
entire solution of system (\ref{11}). The solution constructed in this way
will be radially symmetric. Since the radial solutions of (\ref{11}) are
solutions of the ordinary differential equations system (\ref{77}) it
follows that the radial solutions of (\ref{11}) with $u_{1}\left( 0\right)
=a,$ $u_{2}\left( 0\right) =b$ satisfy:%
\begin{eqnarray}
u_{1}\left( r\right) &=&a+\int_{0}^{r}\frac{1}{t^{N-1}}%
\int_{0}^{t}s^{N-1}p_{1}\left( s\right) f_{1}\left( u_{2}\left( s\right)
\right) dsdt,\text{ }r\geq 0,  \label{eq11} \\
u_{2}\left( r\right) &=&b+\int_{0}^{r}\frac{1}{t^{N-1}}%
\int_{0}^{t}s^{N-1}p_{2}\left( s\right) f_{2}\left( u_{1}\left( s\right)
\right) dsdt,\text{ }r\geq 0.  \label{eq22}
\end{eqnarray}%
Choose $R>0$ so that $r^{2N-2}p_{1}\left( r\right) $ and $%
r^{2N-2}p_{2}\left( r\right) $ are non-decreasing for $r\geq R$. In order to
prove cases i.), ii.), iii.), iv.) and v.) above we intend to establish some
inequalities. Using the same arguments as in (\ref{n8}) and (\ref{nn8}) we
can see that%
\begin{equation}
\left\{ 
\begin{array}{c}
\left[ r^{N-1}\left( u_{1}\left( r\right) \right) ^{\prime }\right] ^{\prime
}\leq r^{N-1}p_{1}\left( r\right) \overline{c}_{1}h_{1}\left(
M_{1}f_{2}\left( u_{1}\left( r\right) \right) \right) \omega _{1}\left(
1+\int_{0}^{r}t^{1-N}\int_{0}^{t}s^{N-1}p_{2}\left( s\right) dsdt\right) ,
\\ 
\left[ r^{N-1}\left( u_{2}\left( r\right) \right) ^{\prime }\right] ^{\prime
}\leq r^{N-1}p_{2}\left( r\right) \overline{c}_{2}h_{2}\left(
M_{2}f_{1}\left( u_{2}\left( r\right) \right) \right) \omega _{2}\left(
1+\int_{0}^{r}t^{1-N}\int_{0}^{t}s^{N-1}p_{1}\left( s\right) dsdt\right) .%
\end{array}%
\right.  \label{n88}
\end{equation}%
Multiplying the first equation in (\ref{n88}) by $r^{N-1}\left( u_{1}\right)
^{\prime }$ and the second by $r^{N-1}\left( u_{2}\right) ^{\prime }$ and
integrating gives 
\begin{equation*}
\left\{ 
\begin{array}{l}
\left[ r^{N-1}\left( u_{1}\left( r\right) \right) ^{\prime }\right] ^{2}\leq %
\left[ R^{N-1}\left( u_{1}\left( R\right) \right) ^{\prime }\right] ^{2} \\ 
+2\int_{R}^{r}z^{2N-2}p_{1}\left( z\right) \overline{c}_{1}\omega _{1}\left(
1+\int_{0}^{z}t^{1-N}\int_{0}^{t}s^{N-1}p_{2}\left( s\right) dsdt\right) 
\frac{d}{dz}\int_{a}^{u_{1}\left( z\right) }h_{1}\left( M_{1}f_{2}\left(
s\right) \right) dsdz, \\ 
\left[ r^{N-1}\left( u_{2}^{k}\left( r\right) \right) ^{\prime }\right]
^{2}\leq \left[ R^{N-1}\left( u_{2}^{k}\left( R\right) \right) ^{\prime }%
\right] ^{2} \\ 
+2\int_{R}^{r}z^{2N-2}p_{1}\left( r\right) \overline{c}_{2}\omega _{2}\left(
1+\int_{0}^{z}t^{1-N}\int_{0}^{t}s^{N-1}p_{1}\left( s\right) dsdt\right) 
\frac{d}{dz}\int_{b}^{u_{2}\left( z\right) }h_{2}\left( M_{2}f_{1}\left(
z\right) \right) dsdz,%
\end{array}%
\right.
\end{equation*}%
for $r\geq R$. We get from the monotonicity of $z^{2N-2}p_{1}\left( z\right) 
$ and $z^{2N-2}p_{2}\left( z\right) $ for $r\geq z\geq R$ that 
\begin{equation*}
\left\{ 
\begin{array}{c}
\left[ r^{N-1}\left( u_{1}\left( r\right) \right) ^{\prime }\right] ^{2}\leq
C_{1}+2\overline{c}_{1}r^{2N-2}p_{1}\left( r\right) \omega _{1}\left(
1+\int_{0}^{r}t^{1-N}\int_{0}^{t}s^{N-1}p_{2}\left( s\right) dsdt\right) 
\overline{H}_{1}\left( u_{1}\left( r\right) \right) , \\ 
\left[ r^{N-1}\left( u_{2}\left( r\right) \right) ^{\prime }\right] ^{2}\leq
C_{2}+2\overline{c}_{2}r^{2N-2}p_{2}\left( r\right) \omega _{2}\left(
1+\int_{0}^{r}t^{1-N}\int_{0}^{t}s^{N-1}p_{1}\left( s\right) dsdt\right) 
\overline{H}_{2}\left( u_{2}\left( r\right) \right) ,%
\end{array}%
\right.
\end{equation*}%
where $C_{1}=\left[ R^{N-1}\left( u_{1}\left( R\right) \right) ^{\prime }%
\right] ^{2}$, $C_{2}=\left[ R^{N-1}\left( u_{2}\left( R\right) \right)
^{\prime }\right] ^{2}$, $\overline{H}_{1}\left( u_{1}\left( r\right)
\right) =\int_{0}^{u_{1}\left( r\right) }h_{1}\left( M_{1}f_{2}\left(
s\right) \right) ds$ and $\overline{H}_{2}\left( u_{2}\left( r\right)
\right) =\int_{0}^{u_{2}\left( r\right) }h_{2}\left( M_{2}f_{1}\left(
s\right) \right) ds$. This implies that%
\begin{equation}
\left\{ 
\begin{array}{c}
\frac{\left( u_{1}\left( r\right) \right) ^{\prime }}{\sqrt{\overline{H}%
_{1}\left( u_{1}\left( r\right) \right) }}\leq \frac{\sqrt{C_{1}}r^{1-N}}{%
\sqrt{\overline{H}_{1}\left( u_{1}\left( r\right) \right) }}+\sqrt{2%
\overline{c}_{1}p_{1}\left( r\right) }\sqrt{\omega _{1}\left(
1+\int_{0}^{r}t^{1-N}\int_{0}^{t}s^{N-1}p_{2}\left( s\right) dsdt\right) },
\\ 
\frac{\left( u_{2}\left( r\right) \right) ^{\prime }}{\sqrt{\overline{H}%
_{2}\left( u_{2}\left( r\right) \right) }}\leq \frac{\sqrt{C_{2}}r^{1-N}}{%
\sqrt{\overline{H}_{2}\left( u_{2}\left( r\right) \right) }}+\sqrt{2%
\overline{c}_{2}p_{2}\left( r\right) }\sqrt{\omega _{2}\left(
1+\int_{0}^{r}t^{1-N}\int_{0}^{t}s^{N-1}p_{1}\left( s\right) dsdt\right) }.%
\end{array}%
\right.  \label{10}
\end{equation}%
In particular, integrating (\ref{10}) from $R$ tor $r$ and using the fact
that 
\begin{eqnarray*}
&&\sqrt{2p_{1}\left( r\right) \overline{c}_{1}\omega _{1}\left(
1+\int_{0}^{r}t^{1-N}\int_{0}^{t}s^{N-1}p_{2}\left( s\right) dsdt\right) } \\
&\leq &r^{1+\varepsilon }p_{1}\left( r\right) \overline{c}_{1}\omega
_{1}\left( 1+\int_{0}^{r}t^{1-N}\int_{0}^{t}s^{N-1}p_{2}\left( s\right)
dsdt\right) +r^{-1-\varepsilon },
\end{eqnarray*}%
lead to%
\begin{equation}
\begin{array}{ll}
\int_{u_{1}\left( R\right) }^{u_{1}\left( r\right) }\left[
\int_{0}^{t}h_{1}\left( M_{1}f_{2}\left( z\right) \right) dz\right] ^{-1/2}dt
& =H_{1}\left( u_{1}\left( r\right) \right) -H_{1}\left( u_{1}\left(
R\right) \right) \\ 
& \leq \sqrt{C_{1}}\int_{R}^{r}t^{1-N}\left( \int_{0}^{u_{1}\left( t\right)
}h_{1}\left( M_{1}f_{2}\left( z\right) \right) dz\right) ^{-1/2}dt \\ 
& +\int_{R}^{r}z^{1+\varepsilon }p_{1}\left( z\right) \overline{c}_{1}\omega
_{1}\left( 1+\int_{0}^{z}\frac{1}{t^{N-1}}\int_{0}^{t}s^{N-1}p_{2}\left(
s\right) dsdt\right) dz \\ 
& +\int_{R}^{r}z^{-1-\varepsilon }dz \\ 
& \leq \frac{\sqrt{C_{1}}\int_{R}^{r}t^{1-N}dt}{\left( \int_{0}^{u_{1}\left(
R\right) }h_{1}\left( M_{1}f_{2}\left( z\right) \right) dz\right) ^{1/2}}+%
\overline{c}_{1}P_{2}\left( r\right) +\frac{1}{\varepsilon R^{\varepsilon }}.%
\end{array}
\label{c1}
\end{equation}%
A special case of this inequality, is originally due to \cite{CDF}. We next
turn to estimating the second sequence. A similar calculation yields 
\begin{equation}
H_{2}\left( u_{2}\left( r\right) \right) -H_{2}\left( u_{2}\left( R\right)
\right) \leq \frac{\sqrt{C_{2}}\int_{R}^{r}t^{1-N}dt}{\left(
\int_{0}^{u_{2}\left( R\right) }h_{2}\left( M_{2}f_{1}\left( z\right)
\right) dz\right) ^{1/2}}+\overline{c}_{2}Q_{2}\left( r\right) +\frac{1}{%
\varepsilon R^{\varepsilon }}.  \label{c2}
\end{equation}%
The inequalities (\ref{c1}) and (\ref{c2}) are needed in proving the
\textquotedblright boundedness\textquotedblright\ of the functions $u_{1}$
and $u_{2}$. Indeed, they can be written as 
\begin{equation}
\left\{ 
\begin{array}{c}
u_{1}\left( r\right) \leq H_{1}^{-1}\left( H_{1}\left( u_{1}\left( R\right)
\right) +\frac{\sqrt{C_{1}}\int_{R}^{r}t^{1-N}dt}{\left(
\int_{0}^{u_{1}\left( R\right) }h_{1}\left( M_{1}f_{2}\left( z\right)
\right) dz\right) ^{1/2}}+\overline{c}_{1}P_{2}\left( r\right) +\frac{1}{%
\varepsilon R^{\varepsilon }}\right) , \\ 
u_{2}\left( r\right) \leq H_{2}^{-1}\left( H_{2}\left( u_{2}\left( R\right)
\right) +\frac{\sqrt{C_{2}}\int_{R}^{r}t^{1-N}dt}{\left(
\int_{0}^{u_{2}\left( R\right) }h_{2}\left( M_{2}f_{1}\left( z\right)
\right) dz\right) ^{1/2}}+\overline{c}_{2}Q_{2}\left( r\right) +\frac{1}{%
\varepsilon R^{\varepsilon }}.\right) .%
\end{array}%
\right.  \label{ineq}
\end{equation}%
Having discussed the~\textquotedblright bounded\textquotedblright\ case, we
now turn to the Cases i.), ii.), iii.), iv.) and v.).

\textbf{Case i.):} When $P_{2}\left( \infty \right) <\infty $ and $%
Q_{2}\left( \infty \right) <\infty $ we find from (\ref{ineq}) that%
\begin{equation*}
\left\{ 
\begin{array}{c}
\lim_{r\rightarrow \infty }u_{1}\left( r\right) <\infty \\ 
\lim_{r\rightarrow \infty }u_{2}\left( r\right) <\infty%
\end{array}%
\right. \text{ for all }r\geq 0
\end{equation*}%
and so $\left( u_{1},u_{2}\right) $ is bounded. \ We next consider:

\textbf{Case ii.):} The case $P_{1}\left( \infty \right) =Q_{1}\left( \infty
\right) =\infty $ is proved in the following: 
\begin{eqnarray}
u_{1}\left( r\right) &=&a+\int_{0}^{r}t^{1-N}\int_{0}^{t}s^{N-1}p_{1}\left(
s\right) f_{1}\left( u_{2}\left( s\right) \right) dsdt  \notag \\
&=&a+\int_{0}^{r}y^{1-N}\int_{0}^{y}t^{N-1}p_{1}\left( t\right) f_{1}\left(
b+\int_{0}^{t}z^{1-N}\int_{0}^{z}s^{N-1}p_{2}\left( s\right) f_{2}\left(
u_{1}\left( s\right) \right) ds)dz\right) dtdy  \notag \\
&\geq &a+\int_{0}^{r}y^{1-N}\int_{0}^{y}t^{N-1}p_{1}\left( t\right)
f_{1}\left( b+f_{2}\left( a\right) \int_{0}^{t}z^{1-N}G_{1}\left( z\right)
dz\right) dtdy  \label{i1} \\
&\geq &\int_{0}^{r}y^{1-N}\int_{0}^{y}t^{N-1}p_{1}\left( t\right)
f_{1}\left( b+f_{2}\left( a\right) \int_{0}^{t}z^{1-N}G_{1}\left( z\right)
dz\right) dtdy  \notag \\
&=&P_{1}\left( r\right) .  \notag
\end{eqnarray}%
Similar arguments show that%
\begin{equation*}
u_{2}\left( r\right) \geq Q_{1}\left( r\right) .
\end{equation*}%
Letting $r\rightarrow \infty $ in (\ref{i1}) and in the above inequality we
conclude that%
\begin{equation*}
\lim_{r\rightarrow \infty }u_{1}\left( r\right) =\lim_{r\rightarrow \infty
}u_{2}\left( r\right) =\infty .
\end{equation*}

\textbf{Case iii.):} In the spirit of Case i.) and Case ii.) above, we have%
\begin{eqnarray*}
\lim_{r\rightarrow \infty }u_{1}\left( r\right) &\leq &H_{1}^{-1}\left( 
\frac{H_{1}\left( u_{1}\left( R\right) \right) \varepsilon R^{\varepsilon }+1%
}{\varepsilon R^{\varepsilon }}+\frac{\sqrt{C_{1}}\int_{R}^{r}t^{1-N}dt}{%
\left( \int_{0}^{u_{1}\left( R\right) }h_{1}\left( M_{1}f_{2}\left( z\right)
\right) zdz\right) ^{1/2}}+\overline{c}_{1}P_{2}\left( r\right) \right) \\
&\leq &H_{1}^{-1}\left( \frac{H_{1}\left( u_{1}\left( R\right) \right)
\varepsilon R^{\varepsilon }+1}{\varepsilon R^{\varepsilon }}+\frac{R^{N-2}%
\sqrt{C_{1}}}{\left( N-2\right) \left( \int_{0}^{u_{1}\left( R\right)
}h_{1}\left( M_{1}f_{2}\left( z\right) \right) dz\right) ^{1/2}}+\overline{c}%
_{1}P_{2}\left( \infty \right) \right) \\
&<&\infty .
\end{eqnarray*}%
Arguing as in \cite{CDF} (see also \cite{CD3}) we have%
\begin{equation*}
\lim_{r\rightarrow \infty }u_{2}\left( r\right) =\infty .
\end{equation*}%
So, if 
\begin{equation*}
P_{2}\left( \infty \right) <\infty \text{ and }Q_{1}\left( \infty \right)
=\infty
\end{equation*}%
we have that 
\begin{equation*}
\lim_{r\rightarrow \infty }u_{1}\left( r\right) <\infty \text{ and }%
\lim_{r\rightarrow \infty }u_{2}\left( r\right) =\infty .
\end{equation*}

\textbf{Case iv.): }By a straightforward modification of the proofs
presented in the Case iii.) the results hold true since any statement about $%
P_{2}\left( \infty \right) $ can be translated into a statement about $%
Q_{2}\left( \infty \right) $.

\textbf{Case v.):} If $\left( u_{1},u_{2}\right) $ is a nonnegative
non-trivial entire large solution of (\ref{11}), then $\left(
u_{1},u_{2}\right) $ satisfy%
\begin{eqnarray}
u_{1}\left( r\right) &\leq &H_{1}^{-1}\left( H_{1}\left( u_{1}\left(
R\right) \right) +\frac{\sqrt{C_{1}}\int_{R}^{r}t^{1-N}dt}{\left(
\int_{0}^{u_{1}\left( R\right) }h_{1}\left( M_{1}f_{2}\left( z\right)
\right) dz\right) ^{1/2}}+\overline{c}_{1}P_{2}\left( r\right) +\frac{1}{%
\varepsilon R^{\varepsilon }}\right) ,  \label{fin1} \\
u_{2}\left( r\right) &\leq &H_{2}^{-1}\left( H_{2}\left( u_{2}\left(
R\right) \right) +\frac{\sqrt{C_{2}}\int_{R}^{r}t^{1-N}dt}{\left(
\int_{0}^{u_{2}\left( R\right) }h_{2}\left( M_{2}f_{1}\left( z\right)
\right) dz\right) ^{1/2}}+\overline{c}_{2}Q_{2}\left( r\right) +\frac{1}{%
\varepsilon R^{\varepsilon }}\right) ,  \label{fin2}
\end{eqnarray}%
where 
\begin{equation*}
C_{1}=\left[ R^{N-1}\left( u_{1}\left( R\right) \right) ^{\prime }\right]
^{2}\text{ and }C_{2}=\left[ R^{N-1}\left( u_{2}\left( R\right) \right)
^{\prime }\right] ^{2}.
\end{equation*}%
\ Next, assuming to the contrary that 
\begin{equation*}
P_{2}\left( \infty \right) <\infty \text{ and }Q_{2}\left( \infty \right)
<\infty ,
\end{equation*}%
then (\ref{13}) yields by taking $r\rightarrow \infty $ in (\ref{fin1}) and (%
\ref{fin2}).

\subparagraph{\textbf{Proof of Theorem \protect\ref{th2} completed}}

It follows from (\ref{pen}) and the conditions of the theorem that%
\begin{equation*}
\begin{array}{l}
H_{1}\left( u_{1}^{k}\left( r\right) \right) \leq \overline{c}%
_{1}P_{3}\left( \infty \right) <\overline{c}_{1}H_{1}\left( \infty \right)
<\infty , \\ 
H_{2}\left( u_{2}^{k}\left( r\right) \right) \leq \overline{c}%
_{2}Q_{3}\left( \infty \right) <\overline{c}_{2}H_{2}\left( \infty \right)
<\infty .%
\end{array}%
\end{equation*}%
On the other hand, since $H_{1}^{-1}$ and $H_{2}^{-1}$ are strictly
increasing on $\left[ 0,\infty \right) $, we find that%
\begin{equation*}
u_{1}^{k}\left( r\right) \leq H_{1}^{-1}\left( \overline{c}%
_{1}^{1/2}P_{3}\left( \infty \right) \right) <\infty \text{ and }%
u_{2}^{k}\left( r\right) \leq H_{2}^{-1}\left( \overline{c}%
_{2}^{1/2}Q_{3}\left( \infty \right) \right) <\infty ,
\end{equation*}%
and then the non-decreasing sequences $\left\{ u_{1}^{k}\left( r\right)
\right\} ^{k\geq 1}$ and $\left\{ u_{2}^{k}\left( r\right) \right\} ^{k\geq
1}$ are bounded above for all $r\geq 0$ and all $k$. Combining these two
facts, we conclude that $\left( u_{1}^{k}\left( r\right) ,u_{2}^{k}\left(
r\right) \right) \rightarrow \left( u_{1}\left( r\right) ,u_{2}\left(
r\right) \right) $ as $k\rightarrow \infty $ and the limit functions $u_{1}$
and $u_{2}$ are positive entire bounded radial solutions of system (\ref{11}%
).\textbf{\ }

\textbf{Case ii.), iii.), iv) and v): }For the proof, we follow the same
steps and arguments as in the proof of Theorem \ref{th1}.

\begin{remark}
Our proof actually gives a slightly stronger result, compared with other
references. Since the conditions of the form 
\begin{equation*}
H_{1}\left( \infty \right) \leq \infty \text{ and }H_{2}\left( \infty
\right) \leq \infty
\end{equation*}%
are firstly introduced here.
\end{remark}

\end{document}